# NEAR-INTEGRATED GARCH SEQUENCES[1]

By István Berkes, Lajos Horváth and Piotr Kokoszka

*Technische Universität Graz and Hungarian Academy of Sciences,*
*University of Utah and Utah State University*

Motivated by regularities observed in time series of returns on speculative assets, we develop an asymptotic theory of GARCH(1, 1) processes $\{y_k\}$ defined by the equations $y_k = \sigma_k \varepsilon_k$, $\sigma_k^2 = \omega + \alpha y_{k-1}^2 + \beta \sigma_{k-1}^2$ for which the sum $\alpha + \beta$ approaches unity as the number of available observations tends to infinity. We call such sequences near-integrated. We show that the asymptotic behavior of near-integrated GARCH(1, 1) processes critically depends on the sign of $\gamma := \alpha + \beta - 1$. We find assumptions under which the solutions exhibit increasing oscillations and show that these oscillations grow approximately like a power function if $\gamma \leq 0$ and exponentially if $\gamma > 0$. We establish an additive representation for the near-integrated GARCH(1, 1) processes which is more convenient to use than the traditional multiplicative Volterra series expansion.

**1. Introduction.** We study the GARCH(1, 1) model of Bollerslev (1986) in which the time series $\{y_k\}$ follows (2.1) and (2.2) of Section 2. The GARCH(1, 1) specification is commonly used as a good approximation for modeling volatility of financial and econometric heteroskedastic time series; see Chapter 15 of Hull (2000) and Chapter 7 of Zivot and Wang (2003) for practical applications and further references. This paper is motivated by the well-established fact that for most time series of returns on speculative assets the sum $\alpha + \beta$ is very close to 1. This observation led Engle and Bollerslev (1986) to introduce the so-called *integrated* GARCH or IGARCH processes which in the setting of (2.1) and (2.2) are defined by the condition $\alpha + \beta = 1$. The term "integrated" is borrowed from the classical theory of linear time series [see, e.g., Brockwell and Davis (1991)], in which a random sequence

---

Received October 2003; revised April 2004.

[1]Supported by NSF Grants INT-02-23262 and DMS-04-13652, NATO Grant PST.EAP.CLG 980599 and Hungarian National Foundation for Scientific Research Grants 43037 and 37886.

*AMS 2000 subject classifications.* 62M10, 91B84.

*Key words and phrases.* Asymptotic distribution, near-integrated GARCH.

$\{x_k\}$ is said to be integrated if the sequence of differences $u_t = x_t - x_{t-1}$ is stationary. If this is the case, the $x_k$ are sums or "integrals" of the $u_k$. The sequence $\{x_k\}$ is then typically not stationary and its realizations resemble a random walk. To see why a GARCH(1,1) process with $\alpha + \beta = 1$ might be called integrated, note that (2.2) may be rewritten as

$$y_k^2 = (\alpha + \beta)y_{k-1}^2 + e_k - \beta e_{k-1}, \qquad e_k := (\varepsilon_k^2 - 1)\sigma_k^2.$$

If there is a strictly stationary solution $\{y_k\}$ to (2.1) and (2.2) and $E\varepsilon_k^2 = 1$, then $\{e_k\}$ is a sequence of martingale differences, so the $y_k^2$ follow an ARMA(1,1) process which is integrated (with $u_k = e_k - \beta e_{k-1}$) if $\alpha + \beta = 1$. The above argument assumed, however, that $\{y_k\}$ was stationary, so $y_k^2$ cannot be nonstationary, unlike in the case of ARMA(1,1) processes with independent innovations $e_k$. In fact, Bougerol and Picard (1992) showed that if $\alpha > 0$ and $\beta > 0$ and $\alpha + \beta = 1$, then the GARCH(1,1) process is strictly stationary.

The objective of the present paper is to investigate the asymptotic behavior of the GARCH(1,1) process when the sum $\alpha + \beta$ is not necessarily equal to 1, but is close to 1, as is the case in applications. This is achieved by allowing both $\alpha$ and $\beta$ to be functions of the sample size $n$; the larger the sample size $n$ is, the closer the sum $\alpha + \beta$ is to unity.

To illustrate some corollaries to our results, set $\gamma = \alpha + \beta - 1$. Note that if $|\gamma| = n^{-q}$ for some $1/2 < q < 1$, then assumptions (2.12) and (2.18) are satisfied. Define $y_n \propto f(n)$ to mean that $y_n/f(n)$ converges weakly to the distribution of the innovations $\varepsilon_k$. Then, assuming that for some $1/2 < q < 1$, $\gamma = -n^{-q}$ in the case $\gamma < 0$ and $\gamma = n^{-q}$ in the case $\gamma > 0$, we conclude that under the assumptions stated in Section 2,

$$\begin{aligned}
&\text{if } \gamma < 0, & y_n &\propto \omega^{1/2}n^{q/2} && \text{(Theorem 2.2)};\\
&\text{if } \gamma = 0, & y_n &\propto \omega^{1/2}n^{1/2} && \text{(Theorem 2.4)};\\
&\text{if } \gamma > 0, & y_n &\propto \omega^{1/2}e^{n/2}n^{q/2} && \text{(Theorem 2.6)}.
\end{aligned}$$

The parameter $\omega$ is defined in (2.2).

Under our assumptions, the sequence $y_n$ is therefore always nonstationary and its "oscillations" increase with the sample size at different rates, depending on the sign of $\gamma$. The case $\gamma = 0$ should be contrasted with the aforementioned result of Bougerol and Picard (1992) who showed that if $\alpha$ and $\beta$ are fixed, then the processes $\{y_n\}$ is stationary but the expected value of its marginal distribution is infinite. Under our assumptions, the process $\{y_n\}$ is mildly explosive but after normalizing by $n^{1/2}$ converges to a random variable with finite fourth moment. Similar interpretation is valid in the case $\gamma < 0$. If $\alpha + \beta > 1$, the "oscillations" increase at an exponential rate.

Our work is somewhat related to Nelson (1990a) who considered approximating diffusion processes by discrete time sequences from the ARCH family. While his main focus was on the exponential GARCH processes, he also considered as an example a GARCH(1,1)-$M$ process which is a random walk type process; we refer to the original work of Nelson (1990a) for the details. In our setting Nelson's conditions (2.28), (2.29) and (2.33) can be stated, respectively, as $n\omega_n \to \omega \geq 0$, $n\gamma_n \to -\theta \in \mathbb{R}$ and $n^{1/2}\alpha_n \to \alpha > 0$. Since his goal was to obtain a continuous time diffusion process as a limit, Nelson needed the parameter $\omega$ also to depend on $n$. In our theory, the parameter $\omega$ plays no role because we investigate what happens when the sum $\alpha + \beta$ approaches 1. Unlike Nelson (1990a), we assume that $n|\gamma_n| \to \infty$ and $n^{1/2}\alpha_n \to 0$. Our results show that assuming different rates of convergence leads to a completely different asymptotic behavior.

In another related work, Kazakevičius, Leipus and Viano (2004) studied the stability of general nonlinear processes related to the ARCH($\infty$) models. Their results apply, however, to situations when the limit is a stationary process. Specialized to GARCH(1, 1) (see their Theorem 4.2 and Remark 2.2), their results imply that $\alpha_n/\gamma_n \to \alpha/\gamma$, where now $\alpha$ and $\gamma < 0$ are the parameters of the limit process, is a sufficient and necessary condition for an appropriately defined convergence. We again refer to the original work of Kazakevičius, Leipus and Viano (2004) for full details which are too complex to be presented here. We study the case when $\gamma_n \to 0$ and, as a result, the finite-dimensional distributions of the GARCH(1, 1) sequences, after appropriate normalization, converge to the finite-dimensional distributions of a process which is no longer GARCH(1, 1).

The term "near-integrated" in the title of the paper is used in analogy with the popular term "near unit root" which, in its simplest form, refers to an AR(1) process $\{x_k\}$ defined by $x_k = \rho x_{k-1} + e_k$, where the $e_k$ are errors. The process $\{x_k\}$ is said to have unit root if $\rho = 1$ and near unit root if $\rho$ tends to 1 with the sample size.

The paper is organized as follows. In Section 2, we state the assumptions and the main results and provide a brief discussion of their significance. The proofs are collected in Sections 3–6.

**2. Main results.** Starting with the initial values $\sigma_0$ and $y_0 = \sigma_0\varepsilon_0$, we define the sequences $\sigma_k, 1 \leq k \leq n$, and $y_k, 1 \leq k \leq n$, by the recursions

$$(2.1) \qquad y_k = \sigma_k\varepsilon_k, \qquad\qquad 1 \leq k \leq n,$$

$$(2.2) \qquad \sigma_k^2 = \omega + \alpha y_{k-1}^2 + \beta\sigma_{k-1}^2, \qquad 1 \leq k \leq n,$$

where

$$(2.3) \qquad\qquad \omega > 0, \qquad \alpha \geq 0, \qquad \beta \geq 0.$$

Throughout this paper we assume that (2.1)–(2.3) hold and that $\varepsilon_0, \varepsilon_1, \varepsilon_2, \ldots, \varepsilon_n$ are independent identically distributed random variables (not necessarily with mean zero) such that the distribution of $\varepsilon_k^2$ is nondegenerate.

Nelson [(1990b), formula (6)] showed that the solutions to (2.1) and (2.2) can be written as

$$(2.4) \quad \sigma_k^2 = \sigma_0^2 \prod_{i=1}^{k} (\beta + \alpha \varepsilon_{k-i}^2) + \omega \left[ 1 + \sum_{j=1}^{k-1} \prod_{i=1}^{j} (\beta + \alpha \varepsilon_{k-i}^2) \right], \qquad 1 \le k \le n.$$

Nelson (1990b) also showed that $\sigma_k^2$ has an a.s. limit, as $k \to \infty$, if and only if

$$(2.5) \qquad E \log(\beta + \alpha \varepsilon_0^2) < 0,$$

assuming that $E|\log(\beta + \alpha \varepsilon_0^2)| < \infty$.

In this paper we assume that $\alpha = \alpha_n$ and $\beta = \beta_n$ and $\alpha \to 0, \beta \to 1$, as $n \to \infty$. In this case, by (2.7) below,

$$\lim_{n \to \infty} E \log(\beta_n + \alpha_n \varepsilon_0^2) = 0.$$

We are interested in the behavior of $\sigma_k^2$ and $y_k$ when

$$(2.6) \qquad k = [nt] \qquad \text{for a fixed } 0 < t \le 1.$$

In fact, we study finite-dimensional distributions, but assumption (2.6) conveys the idea and is extensively appealed to in the proofs.

For reasons explained in the Introduction, we assume that

$$(2.7) \qquad E \varepsilon_0^2 = 1.$$

We will show that the behavior of $\sigma_k^2$ (the behavior of $y_k$ can then be easily derived) critically depends on whether the quantity

$$(2.8) \qquad \gamma = \gamma_n := \alpha + \beta - 1$$

goes to zero from left or right or $\gamma \equiv 0$.

Throughout the paper we use all or some of the following assumptions:

$$(2.9) \qquad E \varepsilon_0^4 \, < \, \infty,$$

$$(2.10) \qquad n^{1/2} \alpha \, \to 0,$$

$$(2.11) \qquad n \alpha \, \to \infty,$$

$$(2.12) \qquad n^{1/2} \gamma \, \to 0.$$

In the following, we will often work with the random variables

$$(2.13) \qquad \xi_j = \varepsilon_j^2 - 1.$$

Our first result establishes an additive representation for $\sigma_k^2$ which should be contrasted with the multiplicative representation (2.4). The representation in Proposition 2.1 is valid for near-integrated GARCH(1, 1) sequences under the stated assumptions. Representation (2.4) is easy to obtain for any GARCH(1, 1) sequences by the repeated application of (2.2). The proofs of the remaining results are based on Proposition 2.1.

PROPOSITION 2.1. *If* (2.6), (2.7) *and* (2.9)–(2.12) *hold, then*

$$\sigma_k^2 = \sigma_0^2 e^{k\gamma}\left(1 + \alpha \sum_{1 \le j \le k} \xi_{k-j} + R_k^{(1)}\right)$$

$$+ \omega\left[1 + \sum_{1 \le j \le k-1} e^{j\gamma}\left(1 + \alpha \sum_{1 \le i \le j} \xi_{k-i} + R_{k,j}^{(2)}\right)(1 + R_{k,j}^{(3)})\right].$$

*The remainder terms satisfy*

$$(2.14) \qquad |R_k^{(1)}| = O_P(k(\alpha^2 + \gamma^2)),$$

$$(2.15) \qquad \max_{1 \le j \le k} |R_{k,j}^{(2)}| = O_P(k\alpha^2),$$

$$(2.16) \qquad \max_{1 \le j \le k} \frac{1}{j \log\log j}|R_{k,j}^{(2)}| = O_P(\alpha^2)$$

*(we set* $\log\log x = 1$ *if* $x < 4$*) and*

$$(2.17) \qquad \max_{1 \le j \le k} \frac{1}{j}|R_{k,j}^{(3)}| = O_P(\alpha^2 + \gamma^2).$$

REMARK 2.1. In relations (2.14)–(2.17) and throughout the paper, the symbols $O_P$ and $o_P$ are meant as $n \to \infty$, although in many places the formulas contain only the letter $k$. However, because of assumption (2.6) in which $t$ is *fixed*, everything in the considered formulas depends in fact on $n$.

We describe now the asymptotic behavior of the vectors $[\sigma_{[nt_m]}^2, m = 1, 2, \ldots, N]$ and $[y_{[nt_m]}, m = 1, 2, \ldots, N]$, where $N$ is a fixed integer and

$$0 < t_1 < t_2 < \cdots < t_N \le 1.$$

We first consider the case $\gamma < 0$. We assume

$$(2.18) \qquad n|\gamma| \to \infty,$$

$$(2.19) \qquad \alpha n^{1/2} \log\log n \to 0,$$

$$(2.20) \qquad \frac{|\gamma|^{3/2}}{\alpha} \to 0,$$

$$(2.21) \qquad E|\varepsilon_0|^{4+\delta} < \infty \qquad \text{for some } \delta > 0.$$

THEOREM 2.1. *Suppose* $\gamma < 0$ *and assumptions* (2.7), (2.11), (2.12) *and* (2.18)–(2.21) *hold. Then the random variables*

$$\frac{\sqrt{2}|\gamma|^3}{\alpha} \frac{1}{\sqrt{E\xi_0^2}}\left[\frac{\sigma_{[nt_m]}^2}{\omega} - \sum_{1 \le j \le [nt_m]-1} e^{j\gamma}\right], \qquad 1 \le m \le N,$$

*are asymptotically independent standard normal.*

THEOREM 2.2. *Under the assumptions of Theorem* 2.1, *the random variables*

$$\left(\frac{|\gamma|}{\omega}\right)^{1/2} y_{[nt_m]}, \qquad 1 \le m \le N,$$

*are asymptotically independent, each with the asymptotic distribution equal to the distribution of* $\varepsilon_0$.

Theorems 2.1 and 2.2 are proved in Section 4.

Next we consider the case when $\gamma = 0$.

THEOREM 2.3. *Suppose* $\gamma = 0$ *and assumptions* (2.7), (2.10), (2.11) *and* (2.21) *hold. Then*

$$\frac{1}{n^{3/2}\alpha}\frac{1}{\sqrt{E\xi_0^2}}\left[\frac{\sigma_{[nt_m]}^2}{\omega} - [nt_m], m = 1, 2, \ldots, N\right] \xrightarrow{d} [\zeta_m, m = 1, 2, \ldots, N],$$

*where* $[\zeta_m, m = 1, 2, \ldots, N]$ *is an* $N$-*variate Gaussian vector with mean zero and covariances*

$$E[\zeta_i\zeta_j] = \tfrac{1}{3}[\min(t_i, t_j)]^3.$$

THEOREM 2.4. *Under the assumptions of Theorem* 2.3, *the random variables*

$$(\omega[nt_m])^{-1/2} y_{[nt_m]}, \qquad 1 \le m \le N,$$

*are asymptotically independent, each with the asymptotic distribution equal to the distribution of* $\varepsilon_0$.

Theorems 2.3 and 2.4 are proved in Section 5.

Finally, we consider the case $\gamma > 0$. In this case, we also need assumptions (2.18) and (2.19), but assumption (2.20) is strengthened to

(2.22) $$\frac{\gamma}{\alpha} = O(1).$$

THEOREM 2.5. *Suppose* $\gamma > 0$ *and* (2.7), (2.9), (2.11), (2.18), (2.19) *and* (2.22) *hold. Then the vector*

$$\left[\frac{1}{(E\xi_0^2)^{1/2}}\frac{\gamma e^{-[nt_m]\gamma}}{\alpha[nt_m]^{1/2}}\left(\frac{\sigma_{[nt_m]}^2}{\omega} - \sum_{1 \le j \le [nt_m]-1} e^{j\gamma}\right), m = 1, 2, \ldots, N\right]$$

*converges in distribution to the vector* $[W(t_m), m = 1, 2, \ldots, N]$, *where* $W(\cdot)$ *is a Wiener process.*

*Estimated parameter values and the implied value of
$\gamma$ based on S&P 500 returns shown in Figure 1*

| $n$ | $\alpha$ | $\beta$ | $\gamma$ |
|------|-----------|-----------|------------|
| 200 | 0.002421982 | 0.993589156 | $-0.003988862$ |
| 300 | 0.0003888147 | 0.9907158442 | $-0.008895341$ |
| 400 | 0.002953501 | 1.002814152 | 0.005767653 |
| 500 | 0.007391279 | 0.992564947 | $-4.3774e-05$ |
| 1000 | 0.020334508 | 0.971043772 | $-0.00862172$ |
| 1500 | 0.03366838 | 0.95160780 | $-0.01472382$ |
| 2000 | 0.065761052 | 0.925748954 | $-0.008489994$ |
| 2500 | 0.057004888 | 0.938578695 | $-0.004416417$ |

THEOREM 2.6. *Under the assumptions of Theorem* 2.5, *the random
variables*

$$(\omega^{-1}\gamma)^{1/2}e^{-[nt_m]\gamma/2}y_{[nt_m]}, \qquad 1 \leq m \leq N,$$

*are asymptotically independent, each with the asymptotic distribution equal
to the distribution of* $\varepsilon_0$.

Theorems 2.5 and 2.6 are proved in Section 6.

We conclude this section with a brief discussion of our results.

Mikosch and Stărică [(2002), Section 3.2] discuss what they call the "IGARCH
effect." Using 9558 daily returns on the S&P 500 index from the mid-
1950s to early 1990s, they estimate the GARCH(1, 1) model on intervals
$I_k = [1, 1500 + 100k]$. Denoting by $\alpha_k$ and $\beta_k$ the estimates based on a re-
alization over the interval $I_k$, they plot the sum $\alpha_k + \beta_k$ against $k$ and
observe that this sum approaches 1 as $k$ increases. This observation and
other arguments lead them to the conclusion that an IGARCH model is not
appropriate for returns data and arises merely from "an accumulation of
nonstationarities." They say that if estimates are based on shorter samples,
then the estimated value of $\gamma$ is not close to zero. We repeated the exper-
iments of Mikosch and Stărică (2002), but used more recent daily returns
on the S&P 500 index starting January 1, 1990. We also worked with the
Dow Jones Industrial Average over the same period, but the results were
very similar, so we focus below only on the S&P 500 index. Figure 1 shows
2500 continuously compounded (logarithmic) daily returns on the S&P 500
index starting January 1, 1990.

Table 1 shows the estimated parameter values and the implied value of
$\gamma$ based on S&P 500 returns over the intervals $\{1, 2, \ldots, n\}$, where time 1
corresponds to January 1, 1990. The estimates were computed using pseudo-
maximum likelihood estimators; that is, the estimates are the values which maximize the likelihood computed under the assumption that the errors $\varepsilon_k$ have standard normal distribution. It is seen from Table 1 that the values of $\alpha$ and $\gamma$ are indeed very close to 0; it is, however, less obvious that the value of $\gamma$ tends to zero with $n$, as claimed by Mikosch and Stărică (2002). Thus, while keeping in mind that the empirical analysis in Mikosch and Stărică (2002) is more thorough than ours, we will satisfy ourselves by saying that the fact the values of $\gamma$ and $\alpha$ are often very close to zero is well documented.

Our results lend themselves to the following interpretation: Rather than postulating that the observed sample $\{y_1, y_2, \ldots, y_N\}$ follows one model, say, $M$, we consider a sequence of models, $M_n$, such that the observations $\{y_1, y_2, \ldots, y_n\}$ follow model $M_n, n \leq N$. This setting appears reasonable, because in practice the equality $\alpha + \beta = 1$ never holds exactly and the case of the exact equality is mathematically very special. On the other hand, if $\alpha$ and $\beta$ are assumed constant and $\alpha + \beta < 1$, mathematically, we deal with the usual stationary case, no matter how close the sum $\alpha + \beta$ is to unity.

Since our objective was to develop a rigorous probabilistic theory for near-integrated GARCH(1, 1) processes and point out some implications of commonly made assumptions, we do not wish to speculate at length on the empirical consequences of our theory; any such attempt must be supported by a more extensive empirical investigation. We merely point out that since $y_n$ follows model $M_n$, the observed sequence $y_n$ should have properties implied by the sequence of models $M_n$. This means that the $y_n$ should exhibit increasing oscillations if the assumptions of the present paper are satisfied. Such oscillations are not apparent for the returns in Figure 1 up to, say, $N = 1000$, yet the estimated values of $\alpha$ and $\gamma$ are very small. This points towards a conclusion that it is not appropriate to use a GARCH(1, 1) model with $\alpha$ close to zero and $\beta$ close to 1 for all samples $\{y_1, y_2, \ldots, y_n\}, n \leq 1000$; such a model is not appropriate at least for some of these samples. Thus our results lend some qualified support to the findings of Mikosch and Stărică (2002), at least as far as modeling returns on indexes is concerned.

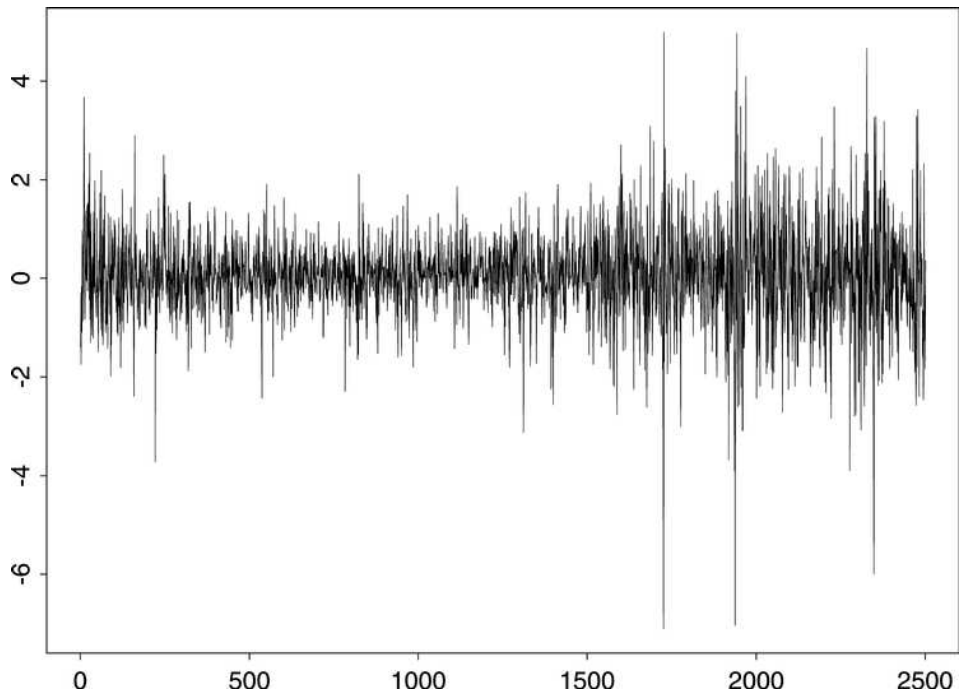

FIG. 1. *Twenty-five hundred continuously compounded daily returns on the S&P 500 index starting January 1, 1990.*

**3. Proof of Proposition 2.1.** In the proofs that follow, we use $a_n \sim b_n$ to denote $\lim_{n\to\infty} a_n/b_n = 1$.

The proof of Proposition 2.1 will use the following elementary result:

LEMMA 3.1. *If* (2.6), (2.7), (2.9), (2.10) *and* (2.12) *hold, then*

$$\max_{1 \leq i < k} |\beta + \alpha \varepsilon_{k-i}^2 - 1| = o_P(1).$$

PROOF. Note that

$$\max_{1 \leq i < k} |\beta + \alpha \varepsilon_{k-i}^2 - 1| \leq |\gamma| + \alpha \max_{1 \leq i < k} |\varepsilon_{k-i}^2 - 1| = |\gamma| + \alpha \max_{1 \leq j \leq k-1} |\varepsilon_j^2 - 1|.$$

Since $\varepsilon_j^2$ are independent and identically distributed, assumption (2.9) and Corollary 3 on page 90 of Chow and Teicher (1988) yield

$$\max_{1 \leq j \leq k-1} |\varepsilon_j^2 - 1| = O(k^{1/2}) \qquad \text{a.s.}$$

Thus the lemma follows from (2.12) and (2.10) combined with (2.6). □

PROOF OF PROPOSITION 2.1. The starting point is formula (2.4). By carefully estimating the remainder terms, we will show in three steps that it can be written as claimed in Proposition 2.1.

*Step* 1. Consider the sequence of events

$$A_n = \left\{ \max_{1 \leq i < k} |\beta + \alpha \varepsilon_{k-i}^2 - 1| \leq \tfrac{1}{2} \right\}.$$

By Lemma 3.1, $\lim_{n\to\infty} P(A_n) = 1$. Since by the Taylor expansion $|\log(1 + x) - x| \leq 2x^2, |x| \leq 1/2$, on the event $A_n$,

$$\left| \sum_{1 \leq i \leq j} \log(\beta + \alpha \varepsilon_{k-i}^2) - \sum_{1 \leq i \leq j} (\gamma + \alpha \xi_{k-i}) \right|$$

$$\leq 2 \sum_{1 \leq i \leq j} (\gamma + \alpha \xi_{k-i})^2 \leq 4\gamma^2 j + 4\alpha^2 \sum_{1 \leq i \leq j} \xi_{k-i}^2.$$

By the law of large numbers and assumption (2.9),

$$\max_{1 \le j \le k} \frac{1}{j} |\xi_{k-1}^2 + \cdots + \xi_{k-j}^2| \stackrel{d}{=} \max_{1 \le j \le k} \frac{1}{j} |\xi_1^2 + \cdots + \xi_j^2| = O_P(1).$$

Therefore

$$\begin{aligned}
(3.1) \qquad \prod_{1 \le i \le j} (\beta + \alpha \varepsilon_{k-i}^2) &= \exp\left[ \sum_{1 \le i \le j} \log(\beta + \alpha \varepsilon_{k-i}^2) \right] \\
&= \exp(\gamma j) \exp\left( \alpha \sum_{1 \le i \le j} \xi_{k-i} \right) \exp(\eta_{k,j}),
\end{aligned}$$

where

$$(3.2) \qquad \max_{1 \le j \le k} \frac{1}{j} |\eta_{k,j}| = O_P(\alpha^2 + \gamma^2).$$

We have thus shown that

$$(3.3) \qquad \prod_{i=1}^{j} (\beta + \alpha \varepsilon_{k-i}^2) = e^{j\gamma} \exp\left( \alpha \sum_{1 \le i \le j} \xi_{k-i} \right) (1 + R_{k,j}^{(3)}),$$

with $R_{k,j}^{(3)}$ satisfying (2.17).

We also note at this point that by (3.2) and assumptions (2.10) and (2.12),

$$(3.4) \qquad \max_{1 \le j \le k} |R_{k,j}^{(3)}| = O_P\left( \max_{1 \le j \le k} |\eta_{k,j}| \right) = O_P(k(\alpha^2 + \gamma^2)) = o_P(1).$$

Formula (3.4) will be used in the verification of (2.14) in step 3.

*Step* 2. In this step we work with the middle factor in the right-hand side of (3.3).

Since $\xi_1, \xi_2, \ldots$ are independent, identically distributed with $E\xi_0^2 < \infty$, by the weak convergence of partial sums to a Brownian motion, we have

$$(3.5) \qquad \max_{1 \le j \le k} \left| \sum_{1 \le i \le j} \xi_{k-i} \right| = O_P(k^{1/2})$$

and therefore by assumption (2.10)

$$(3.6) \qquad \max_{1 \le j \le k} \left| \alpha \sum_{1 \le i \le j} \xi_{k-i} \right| = o_P(1).$$

Let $B_n$ be the event defined by

$$B_n = \left\{ \max_{1 \le j \le k} \left| \alpha \sum_{1 \le i \le j} \xi_{k-i} \right| \le \tfrac{1}{2} \right\}.$$

By (3.6), $\lim_{n\to\infty} P(B_n) = 1$.

Since for $|x| \le 1/2$, $|\exp(x) - (1+x)| \le \sqrt{e}\, x^2/2$, on the event $B_n$

$$(3.7) \qquad \left| \exp\!\left( \alpha \sum_{1 \le i \le j} \xi_{k-i} \right) - \left( 1 + \alpha \sum_{1 \le i \le j} \xi_{k-i} \right) \right| \le \frac{\sqrt{e}}{2} \left( \alpha \sum_{1 \le i \le j} \xi_{k-i} \right)^2.$$

By (3.5),

$$(3.8) \qquad \max_{1 \le j \le k} \left( \alpha \sum_{1 \le i \le j} \xi_{k-i} \right)^2 = O_P(k\alpha^2).$$

By the law of the iterated logarithm

$$(3.9) \qquad \max_{1 \le j \le k} \frac{1}{j \log \log j} \left( \alpha \sum_{1 \le i \le j} \xi_{k-i} \right)^2 = O_P(\alpha^2).$$

We have thus verified that

$$(3.10) \qquad \prod_{i=1}^{j} (\beta + \alpha \varepsilon_{k-i}^2) = e^{j\gamma} \left( 1 + \alpha \sum_{1 \le i \le j} \xi_{k-i} + R_{k,j}^{(2)} \right) (1 + R_{k,j}^{(3)}),$$

with $R_{k,j}^{(2)}$ satisfying (2.15) and (2.16).

*Step* 3. In this step we show that

$$(3.11) \qquad \prod_{i=1}^{k} (\beta + \alpha \varepsilon_{k-i}^2) = e^{k\gamma} \left( 1 + \alpha \sum_{1 \le i \le k} \xi_{k-i} + R_k^{(1)} \right),$$

with $R_k^{(1)}$ satisfying (2.14).

By (3.10) and (3.4),

$$\prod_{1 \le i \le k} (\beta + \alpha \varepsilon_{k-i}^2) = e^{k\gamma} \left( 1 + \alpha \sum_{1 \le i \le k} \xi_{k-i} + O_P(k\alpha^2) \right) (1 + O_P(k(\alpha^2 + \gamma^2))).$$

By the central limit theorem $\sum_{1 \le i \le k} \xi_{k-i} = O_P(k^{1/2})$, so using (2.10) we get

$$\left( 1 + \alpha \sum_{1 \le i \le k} \xi_{k-i} + O_P(k\alpha^2) \right) (1 + O_P(k(\alpha^2 + \gamma^2)))$$

$$= 1 + \alpha \sum_{1 \le i \le k} \xi_{k-i} + O_P(k(\alpha^2 + \gamma^2)),$$

establishing (3.11) with $R_k^{(1)}$ satisfying (2.14).

The proposition now follows from (2.4), (3.10) and (3.11). $\square$

**4. Proof of Theorems 2.1 and 2.2.** Throughout this section we assume that $\gamma = \gamma_n < 0$ and $n|\gamma| \to \infty$.

Using Proposition 2.1, we write for $k = [nt], 0 < t \le 1$,

$$
\begin{aligned}
\sigma_k^2 = {} & \omega + \sigma_0^2 e^{k\gamma} \left( 1 + \alpha \sum_{1 \le j \le k} \xi_{k-j} + R_k^{(1)} \right) \\
& + \omega \sum_{1 \le j \le k-1} e^{j\gamma} \left( 1 + \alpha \sum_{1 \le i \le j} \xi_{k-i} + R_{k,j}^{(2)} \right) R_{k,j}^{(3)} \\
& + \omega \sum_{1 \le j \le k-1} e^{j\gamma} R_{k,j}^{(2)} \\
& + \omega \sum_{1 \le j \le k-1} e^{j\gamma} \left( 1 + \alpha \sum_{1 \le i \le j} \xi_{k-i} \right) \\
=: {} & \omega + \sigma_{k,1}^2 + \sigma_{k,2}^2 + \sigma_{k,3}^2 + \sigma_{k,4}^2.
\end{aligned}
\tag{4.1}
$$

Since $|\gamma|^{3/2}/\alpha \to 0$, the constant term $\omega$ is asymptotically negligible in Theorem 2.1. In Lemma 4.3, we show that the term

$$
\omega^{-1} \sigma_{k,4}^2 = \sum_{1 \le j \le k-1} e^{j\gamma} \left( 1 + \alpha \sum_{1 \le i \le j} \xi_{k-i} \right)
$$

yields the required asymptotic distribution. Since we work with finite-dimensional distributions, we set $k(m) = [nt_m]$ and consider the quantities

$$
\tau_m^* = \sum_{1 \le j \le k(m)-1} e^{j\gamma} \sum_{1 \le i \le j} \xi_{k(m)-i}, \qquad 1 \le m \le N,
\tag{4.2}
$$

which are obtained from $\omega^{-1} \sigma_{[nt_m],4}^2$ after centering and norming.

In Lemmas 4.4, 4.5 and 4.6, respectively, we then show that the terms $\sigma_{k,1}^2, \sigma_{k,2}^2, \sigma_{k,3}^2$ are negligible.

We begin with Lemma 4.1, which is used throughout the proof, and Lemma 4.2, which applies Liapounov's central limit theorem (hence the assumption $E|\varepsilon_0|^{4+\delta} < \infty$) and is then used to prove Lemma 4.3.

Theorem 2.2 follows readily from Theorem 2.1 which together with the relation $\sum_{1 \le j \le k-1} e^{j\gamma} \sim |\gamma|^{-1}$ (see Lemma 4.1) implies

$$
\frac{|\gamma| \sigma_{k(m)}^2}{\omega} - 1 = O_P \left( \frac{\alpha}{|\gamma|^{1/2}} \right) = o_P(1),
$$

because by (2.10) and (2.18), $\alpha/|\gamma|^{1/2} \to 0$. It remains to notice that

$$
\left( \frac{|\gamma|}{\omega} \right)^{1/2} y_{[nt_m]} = \varepsilon_{k(m)} \left( \frac{|\gamma| \sigma_{k(m)}^2}{\omega} \right)^{1/2}
$$

and appeal to the assumption that the $\varepsilon_{k(m)}, 1 \le m \le N$, are independent and identically distributed.

LEMMA 4.1. *For any $\nu \geq 0$,*

$$\sum_{1 \leq j \leq k} j^\nu e^{\gamma j} \sim \frac{1}{|\gamma|^{\nu+1}} \Gamma(\nu + 1).$$

PROOF. Since $x^\nu e^{\gamma x}$ is increasing on $[0, |\gamma|^{-1}\nu]$ and decreasing on $[|\gamma|^{-1}\nu, \infty)$, the sum can be approximated with integrals. Since

$$\int_0^{|\gamma|^{-1}\nu-1} x^\nu e^{\gamma x}\, dx \leq \sum_{1 \leq j \leq |\gamma|^{-1}\nu} j^\nu e^{\gamma j} \leq \int_1^{|\gamma|^{-1}\nu+1} x^\nu e^{\gamma x}\, dx,$$

we conclude that

$$\sum_{1 \leq j \leq |\gamma|^{-1}\nu} j^\nu e^{\gamma j} \sim \frac{1}{|\gamma|^{\nu+1}} \int_0^\nu x^\nu e^{-x}\, dx$$

and a similar argument gives

$$\sum_{|\gamma|^{-1}\nu \leq j \leq k} j^\nu e^{\gamma j} \sim \frac{1}{|\gamma|^{\nu+1}} \int_\nu^\infty x^\nu e^{-x}\, dx. \qquad \square$$

Before formulating our next lemma, we need to introduce new notation. For $0 < t_1 < t_2 < \cdots < t_N < 1$ define $k(m) = [nt_m], 1 \leq m \leq N$, and

(4.3) $$\tau_m = \sum_{1 \leq i \leq k(m)-1} e^{i\gamma} \xi_{k(m)-i}, \qquad 1 \leq m \leq N.$$

LEMMA 4.2. *Suppose that $E|\varepsilon_0|^{4+\delta} < \infty$ for some $\delta > 0$ and that assumptions* (2.7) *and* (2.12) *hold. Then*

$$(2|\gamma|)^{1/2}[\tau_1, \tau_2, \ldots, \tau_N] \xrightarrow{d} (E\xi_0^2)^{1/2}[\eta_1, \eta_2, \ldots, \eta_N],$$

*where $\eta_1, \eta_2, \ldots, \eta_N$ are independent standard normal random variables.*

PROOF. We will use the Cramér–Wold device; see Theorem 29.4 of Billingsley (1995). For any real $\mu_1, \mu_2, \ldots, \mu_N$

$$\mu_1 \tau_1 + \mu_2 \tau_2 + \cdots + \mu_N \tau_N$$
$$= \sum_{1 \leq i \leq k(1)-1} (\mu_1 e^{(k(1)-i)\gamma} + \mu_2 e^{(k(2)-i)\gamma} + \cdots + \mu_N e^{(k(N)-i)\gamma})\xi_i$$
$$\quad + \sum_{k(1) \leq i \leq k(2)-1} (\mu_2 e^{(k(2)-i)\gamma} + \cdots + \mu_N e^{(k(N)-i)\gamma})\xi_i + \cdots$$
$$\quad + \sum_{k(N-1) \leq i \leq k(N)-1} \mu_N e^{(k(N)-i)\gamma}\xi_i$$
$$=: S_1 + S_2 + \cdots + S_N.$$

Observe that

$$
\begin{aligned}
ES_1^2 &= E\xi_0^2 \sum_{1 \le i \le k(1)-1} (\mu_1 e^{(k(1)-i)\gamma} + \mu_2 e^{(k(2)-i)\gamma} + \cdots + \mu_N e^{(k(N)-i)\gamma})^2 \\
&= E\xi_0^2 \sum_{1 \le j \le N} \mu_j^2 \sum_{1 \le i \le k(1)-1} e^{2(k(j)-i)\gamma} \\
&\quad + E\xi_0^2 \sum_{1 \le j \ne l \le N} \mu_j \mu_l \sum_{1 \le i \le k(1)-1} e^{(k(j)+k(l)-2i)\gamma}.
\end{aligned}
$$

By Lemma 4.1,

$$
\tag{4.4}
\sum_{1 \le i \le k(1)-1} e^{2(k(1)-i)\gamma} = \sum_{1 \le i \le k(1)-1} e^{2i\gamma} \sim \frac{1}{2|\gamma|}
$$

and for $2 \le j \le N$,

$$
\tag{4.5}
\sum_{1 \le i \le k(1)-1} e^{2(k(j)-i)\gamma} = e^{2(k(j)-k(1))\gamma} \sum_{1 \le i \le k(1)-1} e^{2(k(1)-i)\gamma} = o\left(\frac{1}{|\gamma|}\right).
$$

Similarly,

$$
\begin{aligned}
\sum_{1 \le i \le k(1)-1} e^{(k(j)+k(l)-2i)\gamma} &= e^{(k(j)-k(1))\gamma} e^{(k(l)-k(1))\gamma} \sum_{1 \le i \le k(1)-1} e^{2(k(1)-i)\gamma} \\
&\sim \frac{1}{2|\gamma|} e^{(k(j)-k(1))\gamma} e^{(k(l)-k(1))\gamma} = o\left(\frac{1}{|\gamma|}\right).
\end{aligned}
$$

Therefore

$$
ES_1^2 = \mu_1^2 E\xi_0^2 \frac{1}{2|\gamma|} + o\left(\frac{1}{|\gamma|}\right).
$$

The same argument applies to $ES_j^2$ for any $1 \le j \le N$, so we see that

$$
\tag{4.6}
\begin{aligned}
E(\mu_1 \tau_1 &+ \mu_2 \tau_2 + \cdots + \mu_N \tau_N)^2 \\
&= (\mu_1^2 + \mu_2^2 + \cdots + \mu_N^2) E\xi_0^2 \frac{1}{2|\gamma|} + o\left(\frac{1}{|\gamma|}\right).
\end{aligned}
$$

Observe also that

$$
\mu_1 \tau_1 + \mu_2 \tau_2 + \cdots + \mu_N \tau_N = \sum_{1 \le i \le k(N)-1} c_i \xi_i
$$

for some $c_i, 1 \le i \le k(N)-1$. It will then follow from (4.6) and Liapounov's central limit theorem [see, e.g., Theorem 27.3 on page 362 of Billingsley (1995)] that

$$
\tag{4.7}
\begin{aligned}
(2|\gamma|)^{1/2} (\mu_1 \tau_1 &+ \mu_2 \tau_2 + \cdots + \mu_N \tau_N) \\
&\xrightarrow{d} (E\xi_0^2)^{1/2} (\mu_1^2 + \mu_2^2 + \cdots + \mu_N^2)^{1/2} \eta,
\end{aligned}
$$

where $\eta$ is a standard normal random variable, provided we have verified that for some $\delta > 0$

$$(4.8) \qquad \frac{(\sum_{1 \leq i \leq k(N)-1} |c_i|^{2+\delta} E|\xi_i|^{2+\delta})^{1/(2+\delta)}}{(\sum_{1 \leq i \leq k(N)-1} c_i^2 E\xi_i^2)^{1/2}} = O\left(\left(\frac{1}{|\gamma|}\right)^{1/(2+\delta)-1/2}\right) = o(1).$$

Relation (4.7) implies the claim in the lemma, so it remains to verify (4.8). By (4.6), the denominator in (4.8) behaves like $(1/|\gamma|)^{1/2}$, so choosing $\delta$ so small that $E|\xi_0|^{2+\delta} < \infty$, it is enough to check that

$$(4.9) \qquad \sum_{1 \leq i \leq k(N)-1} |c_i|^{2+\delta} = O(1/|\gamma|).$$

To illustrate, we will verify (4.9) for the summation range $1 \leq i \leq k(1) - 1$; the remaining ranges of $i$ are dealt with in the same way. For $1 \leq i \leq k(1) - 1$ we get, using Jensen's inequality,

$$\begin{aligned}
|c_i|^{2+\delta} &= |\mu_1 e^{(k(1)-i)\gamma} + \mu_2 e^{(k(2)-i)\gamma} + \cdots + \mu_N e^{(k(N)-i)\gamma}|^{2+\delta} \\
&\leq C_1(N)[|\mu_1|^{2+\delta} e^{(k(1)-i)(2+\delta)\gamma} \\
&\qquad + |\mu_2|^{2+\delta} e^{(k(2)-i)(2+\delta)\gamma} + \cdots + |\mu_N|^{2+\delta} e^{(k(N)-i)(2+\delta)\gamma}].
\end{aligned}$$

Summing the above, we get, using Lemma 4.1 as in (4.4) and (4.5),

$$\sum_{1 \leq i \leq k(N)-1} |c_i|^{2+\delta} \sim C_1(N)|\mu_1|^{2+\delta} \frac{1}{(2+\delta)|\gamma|} + o\left(\frac{1}{|\gamma|}\right) \leq C_2(N)\frac{1}{|\gamma|}. \qquad \square$$

Our next lemma shows that the random variables $\tau_m^*$ defined by (4.2) have the same asymptotic distribution as the $\tau_m, 1 \leq m \leq N$.

LEMMA 4.3. *Suppose that $E|\varepsilon_0|^{4+\delta} < \infty$ for some $\delta > 0$ and that assumptions (2.7) and (2.12) hold. Then*

$$(2|\gamma|^3)^{1/2}[\tau_1^*, \tau_2^*, \ldots, \tau_N^*] \xrightarrow{d} (E\xi_0^2)^{1/2}[\eta_1, \eta_2, \ldots, \eta_N],$$

*where $\eta_1, \eta_2, \ldots, \eta_N$ are independent standard normal random variables.*

PROOF. By Lemma 4.2, it is enough to verify that for each $m$, $|\gamma|^{3/2}\tau_m^* - |\gamma|^{1/2}\tau_m = o_P(1)$. Fix $m$ and to lighten the notation set $k = k(m)$. By rearrangement,

$$\tau_m^* = \sum_{1 \leq i \leq k-1}\left(\sum_{i \leq j \leq k-1} e^{j\gamma}\right)\xi_{k-i}.$$

It is therefore enough to verify that

$$|\gamma|^3 E\left[\sum_{i=1}^{k-1}\left(\sum_{j=i}^{k-1} e^{j\gamma}\right)\xi_{k-i} - \sum_{i=1}^{k-1}|\gamma|^{-1}e^{i\gamma}\xi_{k-i}\right]^2 = o(1).$$

Since

$$(4.10) \qquad \sum_{j=i}^{k-1} e^{j\gamma} = \frac{e^{k\gamma} - e^{i\gamma}}{e^{\gamma} - 1}$$

and $[e^{\gamma} - 1]^{-1} = \gamma^{-1} + O(1)$, we have

$$(4.11) \qquad \sum_{j=i}^{k-1} e^{j\gamma} - |\gamma|^{-1} e^{i\gamma} = e^{k\gamma}(\gamma^{-1} + O(1)) + e^{i\gamma}O(1).$$

Hence, by (4.10), for $1 \le i \le k-1$,

$$(4.12) \qquad \left( \sum_{j=i}^{k-1} e^{j\gamma} - |\gamma|^{-1} e^{i\gamma} \right)^2 = O(|\gamma|^{-1}).$$

By the independence of the $\xi_i$,

$$(4.13) \qquad \begin{aligned} E&\left[ \sum_{i=1}^{k-1} \left( \sum_{j=i}^{k-1} e^{j\gamma} \right) \xi_{k-i} - \sum_{i=1}^{k-1} |\gamma|^{-1} e^{i\gamma} \xi_{k-i} \right]^2 \\ &= E\xi_0^2 \sum_{i=1}^{k-1} \left( \sum_{j=i}^{k-1} e^{j\gamma} - |\gamma|^{-1} e^{i\gamma} \right)^2. \end{aligned}$$

The claim thus follows by (4.12) and $k\gamma^2 \to 0$. $\quad\square$

LEMMA 4.4.   *If* (2.6), (2.7) *and* (2.9)–(2.12) *hold, then*

$$\frac{|\gamma|^{3/2}}{\alpha} \left| e^{k\gamma} \left( 1 + \alpha \sum_{1 \le j \le k} \xi_{k-j} + R_k^{(1)} \right) \right| = o_P(1),$$

*where* $R_k^{(1)}$ *satisfies* (2.14).

PROOF.   Since $k^{-1/2} \sum_{1 \le j \le k} \xi_{k-j}$ is asymptotically normal, conditions (2.10), (2.12) and (2.14) yield

$$\alpha \sum_{1 \le j \le k} \xi_{k-j} + R_k^{(1)} = o_P(1).$$

It remains to observe that by (2.11)

$$\frac{|\gamma|^{3/2}}{\alpha} e^{k\gamma} = o(1)k|\gamma|e^{k\gamma} = o(1)$$

because $\gamma < 0, k|\gamma| \to \infty$. $\quad\square$

LEMMA 4.5. *Suppose* (2.6), (2.7), (2.10) *and* (2.20) *hold. Then*

$$\frac{|\gamma|^{3/2}}{\alpha} \sum_{1 \le j \le k-1} e^{j\gamma} \left| \left( 1 + \alpha \sum_{1 \le i \le j} \xi_{k-i} + R_{k,j}^{(2)} \right) R_{k,j}^{(3)} \right| = o_P(1),$$

*provided* $R_{k,j}^{(2)}$ *satisfies* (2.15) *and* $R_{k,j}^{(3)}$ *satisfies* (2.17).

PROOF. By the weak convergence of partial sums,

$$\max_{1 \le j \le k-1} \left| \sum_{1 \le i \le j} \xi_{k-i} \right| = O_P(k^{1/2}),$$

so by (2.10) and (2.15)

$$(4.14) \qquad \max_{1 \le j \le k-1} \left| \alpha \sum_{1 \le i \le j} \xi_{k-i} + R_{k,j}^{(2)} \right| = o_P(1).$$

By (2.17), (4.14) and Lemma 4.1,

$$\frac{|\gamma|^{3/2}}{\alpha} \sum_{1 \le j \le k-1} e^{j\gamma} \left| \left( 1 + \alpha \sum_{1 \le i \le j} \xi_{k-i} + R_{k,j}^{(2)} \right) R_{k,j}^{(3)} \right|$$

$$= O_P(1) \frac{|\gamma|^{3/2}}{\alpha} (\alpha^2 + \gamma^2) \sum_{1 \le j \le k-1} j e^{j\gamma}$$

$$= O_P(1) \frac{|\gamma|^{3/2}}{\alpha} (\alpha^2 + \gamma^2) \frac{1}{|\gamma|^2}.$$

It remains to observe that

$$\frac{|\gamma|^{3/2}}{\alpha} (\alpha^2 + \gamma^2) \frac{1}{|\gamma|^2} = \frac{\alpha}{|\gamma|^{1/2}} + \frac{|\gamma|^{3/2}}{\alpha}$$

$$= \frac{\alpha n^{1/2}}{(n|\gamma|)^{1/2}} + \frac{|\gamma|^{3/2}}{\alpha} = o(1),$$

by (2.10), (2.18) and (2.20). □

LEMMA 4.6. *Suppose* (2.6) *and* (2.19) *hold. Then*

$$\frac{|\gamma|^{3/2}}{\alpha} \sum_{1 \le j \le k-1} |R_{k,j}^{(2)}| e^{j\gamma} = o_P(1),$$

*provided* $R_{k,j}^{(2)}$ *satisfies* (2.16).

Proof.   Using (2.16), Lemma 4.1, (2.18) and (2.19), we have

$$\frac{|\gamma|^{3/2}}{\alpha} \sum_{1 \le j \le k-1} |R_{k,j}^{(2)}| e^{j\gamma}$$

$$= O_P(1) \frac{|\gamma|^{3/2}}{\alpha} \alpha^2 \sum_{1 \le j \le k-1} j(\log\log j) e^{j\gamma}$$

$$= O_P(1) |\gamma|^{3/2} \alpha (\log\log k) |\gamma|^{-2}$$

$$= o_P(1) \alpha n^{1/2} \log\log n = o_P(1). \qquad \square$$

**5. Proof of Theorems 2.3 and 2.4.**   Throughout this section we assume that $\gamma = 0$.

As in the proof of Theorem 2.1, setting $\gamma = 0$, we have for $k = [nt], 0 < t < 1$,

$$\sigma_k^2 = \omega + \sigma_0^2 \left( 1 + \alpha \sum_{1 \le j \le k} \xi_{k-j} + R_k^{(1)} \right)$$

$$+ \omega \sum_{1 \le j \le k-1} \left( 1 + \alpha \sum_{1 \le i \le j} \xi_{k-i} + R_{k,j}^{(2)} \right) R_{k,j}^{(3)}$$

$$+ \omega \sum_{1 \le j \le k-1} R_{k,j}^{(2)}$$

$$+ \omega \sum_{1 \le j \le k-1} \left( 1 + \alpha \sum_{1 \le i \le j} \xi_{k-i} \right)$$

$$=: \omega + \sigma_{k,1}^2 + \sigma_{k,2}^2 + \sigma_{k,3}^2 + \sigma_{k,4}^2.$$

By (2.11), $n^{3/2}\alpha \to \infty$, so the term $\omega$ is negligible in Theorem 2.3. In Lemmas 5.2–5.4 we show that $n^{-3/2}\alpha^{-1}(\sigma_{k,1}^2 + \sigma_{k,2}^2 + \sigma_{k,3}^2) = o_P(1)$. Therefore

$$(5.1) \qquad \frac{1}{n^{3/2}\alpha}\left( \frac{\sigma_k^2}{\omega} - k \right) = \frac{1}{n^{3/2}} \sum_{1 \le j \le k-1} \sum_{1 \le i \le j} \xi_{k-i} + o_P(1).$$

By (5.1) and Lemma 5.1, the finite-dimensional distributions of the process

$$\left[ \frac{1}{n^{3/2}\alpha} \frac{1}{[E\xi_0^2]^{1/2}} \left( \frac{\sigma_{[nt]}^2}{\omega} - [nt] \right), \ 0 < t < 1 \right]$$

converge to the finite-dimensional distributions of the process $[\int_0^t x\, dW(x), 0 < t < 1]$. Theorem 2.3 thus follows on observing that

$$E\left[ \int_0^t x\, dW(x) \int_0^s x\, dW(x) \right] = \int_0^{\min(t,s)} x^2\, dx = \tfrac{1}{3}[\min(t,s)]^3.$$

Theorem 2.3 and (2.10) imply that $(\omega[nt_m])^{-1}\sigma_{[nt_m]}^2 \xrightarrow{P} 1$ for each $1 \leq m \leq N$, so Theorem 2.4 follows because $(\omega[nt_m])^{-1/2}y_{[nt_m]} = (\omega[nt_m])^{-1/2}\sigma_{[nt_m]}\varepsilon_{[nt_m]}$.

LEMMA 5.1. *If* (2.7) *and* (2.21) *hold, then*

$$n^{-3/2} \sum_{1 \leq j \leq nt-1} \sum_{1 \leq i \leq j} \xi_{k-i} \xrightarrow{d} (E\xi_0^2)^{1/2} \int_0^t x \, dW(x) \qquad in \ D[0,1],$$

*where* $\{W(x), 0 \leq x < \infty\}$ *is a Wiener process.*

PROOF. Note that

$$\sum_{1 \leq j \leq k-1} \sum_{1 \leq i \leq j} \xi_{k-i} = \sum_{1 \leq i \leq k-1} (k-i)\xi_{k-i}$$

$$= \int_0^{k-1} x \, d\left(\sum_{1 \leq j \leq x} \xi_j\right)$$

$$= (k-1) \sum_{1 \leq j \leq k-1} \xi_j - \int_0^{k-1} \left(\sum_{1 \leq j \leq x} \xi_j\right) dx.$$

By (2.7) and the Komlós, Major and Tusnády approximation [see Chapter 1 of Csörgő and Horváth (1993)], there is a Wiener process $W^*$ such that

$$\sum_{1 \leq j \leq x} \xi_j - (E\xi_0^2)^{1/2}W^*(x) = o(x^{1/(2+\delta)}) \qquad \text{a.s.}$$

Hence

(5.2) $$n^{-3/2} \sup_{0 \leq t \leq 1} \left| \int_0^{nt-1} \left( \sum_{1 \leq j \leq x} \xi_j - (E\xi_0^2)^{1/2}W^*(x) \right) dx \right|$$
$$= O_P(1)n^{-3/2}n^{1+1/(2+\delta)} = o_P(1)$$

and

(5.3) $$n^{-3/2} \sup_{0 \leq t \leq 1} \left| (nt-1) \left( \sum_{1 \leq j \leq nt-1} \xi_j - (E\xi_0^2)^{1/2}W^*(nt-1) \right) \right| = o_P(1).$$

In (5.2) and (5.3) we can clearly replace $nt-1$ by $nt$. To complete the proof it is therefore enough to notice that, by the scale transformation of the Wiener process,

$$\left\{ n^{-3/2}\left( ntW^*(nt) - \int_0^{nt} W^*(x) \, dx \right), \ 0 \leq t \leq 1 \right\}$$

$$\overset{d}{=} \left\{ tW(t) - \int_0^t W(x) \, dx, 0 \leq t \leq 1 \right\}$$

and use integration by parts. □

LEMMA 5.2.   *If* (2.6), (2.7), (2.9)–(2.11) *hold, then*

$$n^{-3/2}\alpha^{-1}\left|1 + \alpha \sum_{1 \leq j \leq k} \xi_{k-j} + R_k^{(1)}\right| = o_P(1),$$

*where* $R_k^{(1)}$ *satisfies* (2.14).

PROOF.   By (2.7) and (2.9), $\sum_{1 \leq j \leq k} \xi_{k-j} = O_P(k^{1/2})$, so the lemma follows immediately from (2.6), (2.10) and (2.11).   □

LEMMA 5.3.   *If* (2.6), (2.7) *and* (2.10) *hold, then*

$$n^{-3/2}\alpha^{-1}\left|\sum_{1 \leq j \leq k-1}\left(1 + \alpha \sum_{1 \leq i \leq j} \xi_{k-i} + R_{k,j}^{(2)}\right)R_{k,j}^{(3)}\right| = o_P(1),$$

*where* $R_{k,j}^{(2)}$ *satisfies* (2.15) *and* $R_{k,j}^{(3)}$ *satisfies* (2.17).

PROOF.   By (4.14), it suffices to verify that $n^{-3/2}\alpha^{-1}\sum_{1 \leq j \leq k-1}|R_{k,j}^{(3)}| = o_P(1)$ which follows immediately from (2.17), (2.10), (2.12) and (2.20). Recall that to establish (4.14) we needed assumptions (2.6), (2.7), (2.10) and relation (2.15).   □

LEMMA 5.4.   *If* (2.6) *and* (2.10) *hold, then*

$$n^{-3/2}\alpha^{-1}\left|\sum_{1 \leq j \leq k-1} R_{k,j}^{(2)}\right| = o_P(1),$$

*where* $R_{k,j}^{(2)}$ *satisfies* (2.15).

PROOF.   Follows immediately from (2.15), (2.6) and (2.10).   □

**6. Proof of Theorems 2.5 and 2.6.**   In this section we assume that $\gamma > 0$ and $n\gamma \to \infty$. Note that the assumptions of Theorem 2.5 imply that all assumptions (2.9)–(2.12) hold. In particular, assumption (2.12) is implied by (2.19) and (2.22).

We again use decomposition (4.1). Lemmas 6.1–6.3 below imply that

$$\frac{\gamma e^{-k\gamma}}{\alpha k^{1/2}}(\omega + \sigma_{k,1}^2 + \sigma_{k,2}^2 + \sigma_{k,3}^2) = o_P(1).$$

Therefore

$$\frac{\gamma e^{-k\gamma}}{\alpha k^{1/2}}\left[\frac{\sigma_k^2}{\omega} - \sum_{1 \leq j \leq k-1} e^{j\gamma}\right] = \frac{\gamma e^{-k\gamma}}{k^{1/2}}\sum_{1 \leq j \leq k-1} e^{j\gamma}\sum_{1 \leq i \leq j} \xi_{k-i} + o_P(1).$$

By Lemma 6.4, the last relation implies

$$\frac{\gamma e^{-k\gamma}}{\alpha k^{1/2}}\left[\frac{\sigma_k^2}{\omega}-\sum_{1\le j\le k-1}e^{j\gamma}\right]=k^{-1/2}\sum_{1\le i\le k-1}\xi_i+o_P(1),$$

so by Donsker's theorem [see, e.g., Theorem 14.1 in Billingsley (1999)], we conclude that the finite-dimensional distributions of the process

$$\left[\frac{1}{(E\xi_0^2)^{1/2}}\frac{\gamma e^{-[nt]\gamma}}{\alpha[nt]^{1/2}}\left[\frac{\sigma_{[nt]}^2}{\omega}-\sum_{1\le j\le [nt]-1}e^{j\gamma}\right],0<t<1\right]$$

converge to the finite-dimensional distributions of the Wiener process $[W(t),0<t<1]$. This completes the proof of Theorem 2.5.

Theorem 2.6 now readily follows. For any $0<t<1$,

$$(\omega^{-1}\gamma)^{1/2}e^{-[nt]\gamma/2}y_{[nt]}=(\omega^{-1}\gamma e^{-[nt]\gamma}\sigma_{[nt]}^2)^{1/2}\varepsilon_{[nt]}. \tag{6.1}$$

Since $\alpha[nt]^{1/2}\to 0$, it follows from Theorem 2.5 that

$$\gamma e^{-[nt]\gamma}\left(\omega^{-1}\sigma_{[nt]}^2-\sum_{1\le j\le k-1}e^{j\gamma}\right)=o_P(1). \tag{6.2}$$

Direct verification shows that

$$\gamma e^{-[nt]\gamma}\left(\sum_{1\le j\le k-1}e^{j\gamma}-\gamma^{-1}e^{[nt]\gamma}\right)=o(1). \tag{6.3}$$

Relations (6.2) and (6.3) yield $\omega^{-1}\gamma e^{-[nt]\gamma}\sigma_{[nt]}^2\xrightarrow{P}1$, which combined with (6.1) concludes the proof of Theorem 2.6.

LEMMA 6.1. *If* (2.6), (2.7), (2.9)–(2.12) *and* (2.22) *hold, then*

$$\frac{\gamma e^{-k\gamma}}{\alpha k^{1/2}}\left|e^{k\gamma}\left(1+\alpha\sum_{1\le j\le k}\xi_{k-j}+R_k^{(1)}\right)\right|=o_P(1),$$

*where $R_k^{(1)}$ satisfies* (2.14).

PROOF. Relation (2.14) and the asymptotic normality of $k^{-1/2}\sum_{1\le j\le k}\xi_{k-j}$ yield

$$\frac{\gamma e^{-k\gamma}}{\alpha k^{1/2}}\left|e^{k\gamma}\left(1+\alpha\sum_{1\le j\le k}\xi_{k-j}+R_k^{(1)}\right)\right|$$

$$=O_P(1)\frac{\gamma}{\alpha k^{1/2}}(1+\alpha k^{1/2}+k(\alpha^2+\gamma^2))$$

$$=O_P(1)\left(\frac{\gamma}{\alpha}k^{-1/2}+\gamma+\gamma(k^{1/2}\alpha)+\frac{(k^{1/2}\gamma)^3}{k\alpha}\right)=o_P(1),$$

by (2.10)–(2.12) and (2.22). □

LEMMA 6.2. *If* (2.6), (2.7), (2.9), (2.10), (2.12) *and* (2.22) *hold, then*

$$\frac{\gamma e^{-k\gamma}}{\alpha k^{1/2}} \left| \sum_{1 \le j \le k-1} e^{j\gamma} \left( 1 + \alpha \sum_{1 \le i \le j} \xi_{k-i} + R_{k,j}^{(2)} \right) R_{k,j}^{(3)} \right| = o_P(1),$$

*where* $R_{k,j}^{(2)}$ *satisfies* (2.15) *and* $R_{k,j}^{(3)}$ *satisfies* (2.17).

PROOF.   By (2.17) and (4.14), it suffices to show that the quantity

$$V_k = \frac{\gamma e^{-k\gamma}}{\alpha k^{1/2}} \sum_{1 \le j \le k-1} e^{j\gamma} j(\alpha^2 + \gamma^2)$$

tends to zero. Since for $\gamma > 0$,

$$(6.4) \qquad\qquad \frac{e^{k\gamma} - e^\gamma}{e^\gamma - 1} < \frac{e^{k\gamma}}{\gamma},$$

we have

$$V_k < k^{1/2}\alpha + (k^{1/2}\gamma)\frac{\gamma}{\alpha} \to 0,$$

by (2.10), (2.12) and (2.22). □

LEMMA 6.3. *If* (2.6), (2.9) *and* (2.19) *hold, then*

$$\frac{\gamma e^{-k\gamma}}{\alpha k^{1/2}} \sum_{1 \le j \le k-1} e^{j\gamma} |R_{k,j}^{(2)}| = o_P(1),$$

*where* $R_{k,j}^{(2)}$ *satisfies* (2.16).

PROOF.   By (2.16), we need to show that the quantity

$$U_k = \frac{\gamma e^{-k\gamma}}{\alpha k^{1/2}} \sum_{1 \le j \le k-1} e^{j\gamma} \alpha^2 j \log \log j$$

tends to zero. By (6.4),

$$U_k \le \frac{\gamma e^{-k\gamma}}{\alpha k^{1/2}} \alpha^2 k \log \log k \frac{e^{k\gamma}}{\gamma} = \alpha k^{1/2} \log \log k = o(1),$$

by (2.19). □

LEMMA 6.4. *If* (2.6), (2.7) *and* (2.9) *hold, then*

$$\frac{\gamma^2}{k} e^{-2k\gamma} E\left( \sum_{1 \le j \le k-1} e^{j\gamma} \sum_{1 \le i \le j} \xi_{k-i} - \frac{e^{k\gamma}}{\gamma} \sum_{1 \le i \le k-1} \xi_i \right)^2 \to 0.$$

PROOF. Note that

$$\sum_{1 \le j \le k-1} e^{j\gamma} \sum_{1 \le i \le j} \xi_{k-i} = \sum_{1 \le i \le k-1} \left( \sum_{i \le j \le k-1} e^{j\gamma} \right) \xi_{k-i}$$

and

$$\sum_{1 \le i \le k-1} \xi_i = \sum_{1 \le i \le k-1} \xi_{k-i}.$$

Therefore,

$$E \left( \sum_{1 \le j \le k-1} e^{j\gamma} \sum_{1 \le i \le j} \xi_{k-i} - \sum_{1 \le i \le k-1} \frac{e^{k\gamma}}{\gamma} \xi_i \right)^2$$

$$= E\xi_0^2 \sum_{1 \le i \le k-1} \left( \frac{e^{k\gamma} - e^{i\gamma}}{e^\gamma - 1} - \frac{e^{k\gamma}}{\gamma} \right)^2.$$

By the Taylor expansion,

$$\left| \frac{e^{k\gamma} - e^{i\gamma}}{e^\gamma - 1} - \frac{e^{k\gamma}}{\gamma} \right| \le C_1 \left( \frac{e^{i\gamma}}{\gamma} + e^{k\gamma} \right).$$

Therefore,

$$\sum_{1 \le i \le k-1} \left( \frac{e^{k\gamma} - e^{i\gamma}}{e^\gamma - 1} - \frac{e^{k\gamma}}{\gamma} \right)^2$$

$$\le 2C_1^2 \left[ \sum_{1 \le i \le k-1} \frac{e^{2i\gamma}}{\gamma^2} + k e^{2k\gamma} \right]$$

$$= O(1) \left[ \frac{1}{\gamma^2} \frac{e^{2k\gamma} - e^{2\gamma}}{e^{2\gamma} - 1} + k e^{2k\gamma} \right]$$

$$= O(1) \left[ \frac{1}{\gamma^3} e^{2k\gamma} + k e^{2k\gamma} \right].$$

The claim follows on observing that

$$\frac{\gamma^2}{k} e^{-2k\gamma} \frac{1}{\gamma^3} e^{2k\gamma} = \frac{1}{\gamma k} \to 0$$

and

$$\frac{\gamma^2}{k} e^{-2k\gamma} k e^{2k\gamma} = \gamma^2 \to 0. \qquad \square$$

I. BERKES
INSTITUT FÜR STATISTIK
TECHNISCHE UNIVERSITÄT GRAZ
STEYRERGASSE 17/IV
8010 GRAZ
AUSTRIA
AND
A. RÉNYI INSTITUTE OF MATHEMATICS
HUNGARIAN ACADEMY OF SCIENCES
P.O. BOX 127
H-1364 BUDAPEST
HUNGARY
E-MAIL: berkes@stat.tu-graz.ac.at

L. HORVÁTH
DEPARTMENT OF MATHEMATICS
UNIVERSITY OF UTAH
155 SOUTH 1440 EAST
SALT LAKE CITY, UTAH 84112
USA
E-MAIL: horvath@math.utah.edu

P. KOKOSZKA
DEPARTMENT OF MATHEMATICS AND STATISTICS
UTAH STATE UNIVERSITY
3900 OLD MAIN HILL
LOGAN, UTAH 84322
USA
E-MAIL: piotr@math.usu.edu